\documentclass[a4paper,12pt]{amsart}

\usepackage{amsmath,graphics}
\usepackage{amssymb}
\usepackage{amsfonts}
\usepackage{latexsym}
\usepackage{eucal}
\usepackage[dvips]{graphicx}
\usepackage{color}
\usepackage{enumerate}
\usepackage{mathrsfs}
\newtheorem{thm}{Theorem}[section]

\newtheorem{propo}[thm]{Proposition}
\newtheorem{lem}[thm]{Lemma}

\renewcommand{\Re}{{\rm Re}}
\renewcommand{\Im}{{\rm Im}}
\newcommand{\R}{\mathbb{R}}
\newcommand{\C}{\mathbb{C}}
\newcommand{\Z}{\mathbb{Z}}
\newcommand{\N}{\mathbb{N}}

\newcommand{\D}{ \mathcal{D}}

\newcommand{\hh}{\mathbb{H}^{2}}

\newcommand{\G}{{\bf G}}

\newcommand{\lt}{{\mathcal L}}

\newcommand{\half}{{\textstyle{\frac{1}{2}}}}
\newcommand{\rr}{{ \varrho}}

\begin{document}
\bibliographystyle{plain}
\title[Abelian covers and infinite volume hyperbolic surfaces]{Spectral gaps and abelian covers of convex co-compact surfaces}
\keywords{Fuchsian groups, Hyperbolic surfaces, Laplacian, Resonances, Selberg zeta functions, Representation theory}

\author{Fr\'ed\'eric Naud}
\address{Laboratoire de Math\'ematiques d'Avignon \\
Campus Jean-Henri Fabre, 301 rue Baruch de Spinoza\\
84916 Avignon Cedex 9, France. }
\email{frederic.naud@univ-avignon.fr}

 \maketitle
\begin{abstract} Given $X=\Gamma \backslash \hh$ a convex co-compact hyperbolic surface, we investigate the resonance spectrum $\mathcal{R}_j$ of the laplacian $\Delta_j$ on
large finite abelian covers $X_j=\Gamma_j\backslash \hh\rightarrow X$, where $\Gamma_j$ is a normal subgroup of $\Gamma$ such that $\Gamma/\Gamma_j:=\G_j$ is a finite abelian group.
We show that there exist a {\it uniform} $\epsilon>0$ such that for all $j$, $\Delta_j$ has only finitely many resonances in the strip $\{ \delta-\epsilon\leq \Re(s)\leq \delta \}$ which are {\it all real } and satisfy a {\it Weyl law} as $j\rightarrow \infty$,
$$\# \mathcal{R}_j\cap\{ \delta-\epsilon\leq \Re(s)\leq \delta \}\sim C\vert \G_j\vert,$$
for some $C>0$. This result is an abelian analog of \cite{OhWinter}, and strengthens a previous local result of \cite{JNS}.
 \end{abstract}
 \tableofcontents
 \section{Introduction and results}
 The asymptotic spectrum of large covers of compact Riemannian manifolds is a rich subject with a history of results that reach far beyond the specialized topic
 of spectral geometry with deep interactions with graph theory, representation theory and number theory. 
 Let us be more precise. Denote by $M$ a closed smooth Riemannian manifold, and 
 let $$\lambda_0(M)=0<\lambda_1(M)\leq \lambda_2(M)\leq\ldots\lambda_n(M)\leq$$ 
 be the sequence of eigenvalues of the positive Laplacian on $M$. Let $\pi_1(M)$ be the fundamental group of $M$ and consider a sequence $(N_j)_{j\leq 1}$ of 
 finite index normal subgroups of $\pi_1(M)$, and denote by $(M_j)$ a corresponding sequence of finite sheeted riemannian covers of $M$, with $\pi_1(M_j)=N_j$.
 Let $S\subset \pi_1(M)$ be a symmetric system of generators of $\pi_1(M)$ and consider the sequence $\mathcal{G}_j(S)$ of {\it Cayley graphs} of the galois group $\G_j=\pi_1(M)/N_j$ with respect to $S$. Then the following are equivalent:
 \begin{enumerate}
 \item There exists $\epsilon_1>0$ such that for all $j$, $\lambda_1(M_j)\geq \epsilon_1$.
 \item The sequence of graphs $\mathcal{G}_j(S)$ is a family of expanders, i.e. there exists $\epsilon_2>0$ such that the combinatorial laplacian on
 $\mathcal{G}_j(S)$ has a uniform spectral gap: for all $j$, $\lambda_1(\mathcal{G}_j(S))\geq \epsilon_2$.
 \end{enumerate}
 This result, that we should attribute to Brooks \cite{Brooks}, combines previous important results by Cheeger \cite{Cheeger}, Buser \cite{Buser}, Alon-Milman \cite{AlonMilman}. An interesting corollary is that if the galois groups $\G_j$ are abelian, because Cayley graphs of finite abelian groups cannot be expanders, then up to a sequence extraction we have 
 $$\lim_{j\rightarrow \infty}\lambda_1(M_j)=0.$$
 This fact that had already been observed by Randol \cite{Randol} for compact hyperbolic surfaces and Selberg \cite{Selberg} for some abelian covers of the modular surface.
 We also mention the recent book of Bergeron \cite{Bergeron} where a simple geometric proof, based on Buser's inequality for compact surfaces, is provided.
 
 In this paper we will address a similar problem for a class of {\it infinite volume} hyperbolic surfaces called convex co-compact. The $L^2$-eigenvalues of the Laplacian will be replaced by {\it resonances}. Let us be more specific. Let $\hh$ be the hyperbolic plane endowed with its metric of constant curvature $-1$.
Let $\Gamma$ be a geometrically finite Fuchsian group of isometries acting on $\hh$. This means
that $\Gamma$ admits a finite sided polygonal fundamental domain in $\hh$. We will require that $\Gamma$ has no {\it elliptic} elements different from the identity and that the quotient $\Gamma \backslash \hh$ is of {\it infinite hyperbolic area}. We assume in addition in this
paper that $\Gamma$ has no parabolic elements (no cusps).
Under these assumptions, the quotient space 
$X=\Gamma \backslash \hh$ is a Riemann surface (called convex co-compact) whose {\it ends geometry} is as follows.
The surface $X$ can be decomposed into a compact surface $N$ with geodesic boundary, called the Nielsen region, on which infinite area ends $F_i$ are glued: the funnels. A funnel $F_i$ is a half cylinder isometric to 
$$F_i=(\R /l_i \Z)_\theta \times (\R^+)_t,$$ where $l_i>0$, with the warped metric 
$$ds^2=dt^2+\cosh^2(t)d\theta^2.$$
The limit set $\Lambda(\Gamma)$ of the group $\Gamma$ is defined as 
$$\Lambda(\Gamma):=\overline{\Gamma.z}\cap \partial \hh,$$
where $z\in \hh$ is a given point and $\Gamma.z$ is the orbit under the action of $\Gamma$ which accumulates
on the boundary $\partial \hh$. The limit set $\Lambda(\Gamma)$ does not depend on the choice of the initial point $z$ and its Hausdorff dimension $\delta(\Gamma)$
is the critical exponent of Poincar\'e series \cite{Patterson}. 

\bigskip
The spectrum of $\Delta_X$ on $L^2(X)$ has been described in the works of
Lax and Phillips and Patterson in \cite{LP1,Patterson} as follows: 
\begin{itemize}
\item The half line $[1/4, +\infty)$ is the absolutely continuous spectrum.
\item There are  no $L^2$ embedded eigenvalues inside $[1/4,+\infty)$.
\item The pure point spectrum is empty if $\delta\leq \half$, and finite and starting at $\delta(1-\delta)$ if $\delta>\half$.
\end{itemize}

As a consequence, the resolvent 
$$R(s):=(\Delta_X-s(1-s) )^{-1}:L^2(X)\rightarrow L^2(X)$$
is a holomorphic family for $\Re(s) >\half$, except at a finite number of possible poles related to the eigenvalues. From the work of Mazzeo-Melrose \cite{MM}, it can be meromorphically continued (to all $\C$)  from $C_0^\infty(X)\rightarrow C^\infty(X)$, and poles are called {\it resonances}. We denote
in the sequel by $\mathcal{R}_X$ the set of resonances, written with multiplicities.

From the convergence of Poincar\'e series and the hypergeometric representation of the Schwartz kernel of the resolvent one can deduce that we always have (see for example in
\cite{GuiNaud})
$$\mathcal{R}_X\subset \{ \Re(s)\leq \delta\},$$
and resonances for which the real part is close to $\delta$ are called {\it non trivial sharp resonances}. They correspond to metastable states with the longest lifetime.
There is already a large history of works studying sharp resonances in the context of hyperbolic geometry, we refer the reader for example to \cite{Zworski} for a survey.
Among several results and conjectures related to density, spectral gaps of the resonance spectrum, perhaps the simplest and still widely open problem is to describe
accurately resonances that lie in a small strip close to the vertical line $\Re(s)=\delta$. In particular, no analog of the inequalities of Cheeger and Buser is known for resonances.
In \cite{Naud2}, the author has shown that there exists a small $\epsilon(\Gamma)>0$ such that
$$\mathcal{R}_X\cap\{ \Re(s)\geq \delta-\epsilon \}=\{\delta\},$$
whenever $X=\Gamma \backslash \hh$ with $\Gamma$ non-elementary, i.e. that there exists a {\it spectral gap} beyond the leading resonance $\delta$. However the size of $\epsilon(\Gamma)$ is barely explicit and its relationship with the geometry of the surface is unknown. In this paper we will prove a precise result that goes into the direction of the results of Brooks \cite{Brooks}, despite the fact that {\it no Cheeger inequality} is known.
Let $\Gamma$ be a convex co-compact group as above, and consider a sequence $(\Gamma_j)$ of finite index normal subgroups of $\Gamma$ such that
$$\G_j:=\Gamma/\Gamma_j$$
is a finite Abelian group. More precisely, $\G_j$ has the following structure
$$\G_j=\Z/N_1^{(j)}\Z\times \ldots\times  \Z/N_k^{(j)} \Z,$$
where $N_1(j),N_2(j),\ldots N_k(j)$ are integers such that
$$\lim_{j\rightarrow +\infty}\min_{\ell=1,\ldots,k} N_\ell^{(j)}=+\infty,$$
see $\S 2$ for more details on the construction we use. 
Very much like in the compact picture, we have a sequence of Galois covers $X_j:=\Gamma_j\backslash \hh\rightarrow X$, and we would like to understand
the behaviour of sharp resonances as $\vert \G_j \vert\rightarrow +\infty$. Notice that since $\Gamma_j$ are finite index subgroups of $\Gamma$ we have 
$\delta(\Gamma_j)=\delta(\Gamma)$ for all $j$, so the "base" resonance remains unchanged.

Our main result is the following.

\begin{thm}
\label{theo1}
Assume that $\Gamma$ is non elementary and consider a sequence of abelian covers as above.
\begin{itemize}
\item Then there exists $\epsilon_0(\Gamma)>0$ such that for all $j\in \N$, 
$$\mathcal{R}_{X_j}\cap \{ \Re(s)\geq \delta-\epsilon_0\}$$
consists of finitely many resonances included in the segment $[\delta-\epsilon_0,\delta]$. 
\item Moreover, up to a sequence extraction, we have weak convergence in $C^0([\delta-\epsilon_0,\delta])^*$ of the spectral measures:
$$\lim_{j\rightarrow +\infty} \frac{1}{\vert \G_j\vert} \sum_{\lambda \in \mathcal{R}_j\cap [\delta-\epsilon_0,\delta]} {\mathcal D }_\lambda=\mu,$$
where $\mu$ is an absolutely continuous finite measure supported on $[\delta-\epsilon_0,\delta]$, and $\mathcal{D}_\lambda$ is the Dirac measure at $\lambda$.
\item In addition, if $\lambda \in \mathcal{R}_{X}$, then for all $\varepsilon_0>0$ small enough, one can find $C_0>0$ such that as $j\rightarrow +\infty$,
$$C_0^{-1}\vert \G_j\vert \leq \#\mathcal{R}_{X_j}\cap D(\lambda,\varepsilon_0)\leq C_0\vert \G_j\vert.$$
\end{itemize}
\end{thm}
The limit measure $\mu$ depends on the sequence of covers, see $\S 3$ for details.
This theorem is the perfect Abelian analog of the main result of Oh-Winter \cite{OhWinter} on congruence subgroups of convex co-compact subgroups of $SL_2(\Z)$.
 We point out that the second part of the theorem was essentially obtained in \cite{JNS}, but in a small neighbourhood of $\delta$. The above theorem gives now a complete picture of sharp resonances in a tiny strip near $\Re(s)=\delta$ for
sequences of Abelian covers. The third part describes the behaviour deeper into the spectrum in the vicinity of each resonance in
$\mathcal{R}_X$. Note that it does not say anything about resonances in the cover away from the preexisting resonances of $X$.

This result is not only interesting for itself as it provides new insights into the structure of resonance spectrum, but it has virtually many applications.
For example, more precise local wave asymptotics as in \cite{GuiNaud}, uniform counting asymptotics for the family of groups $\Gamma_j$ as in \cite{BGS,OhWinter}, and also refined prime orbit theorems with precise dependence on $j$, using zeta function techniques. Another byproduct of our analysis is a new result
for counting closed geodesics in homology classes which extends a previous work of McGowan and Perry \cite{McP} without assumption on $\delta$, see section $\S 2.3$.

The plan of the proof is as follows. We first discuss the general structure of abelian covers and show, using a factorization formula for the selberg zeta function involving
$L$-factors related to characters of $\G_j$, one can reduce the problem to the study of zeros of a multivariate entire function $\mathcal{Z}_\Gamma(s,\theta)$, $s\in \C$, $\theta \in \R^r$ related to the homology $H^1(X,\Z)$, which is done in $\S 2$. The goal is then to show that one can find a small $\epsilon>0$ such that for all $\theta$ and all 
$s \in \{ \delta-\epsilon \leq \Re(s)\leq \delta\}$, the holomorphic function $s\mapsto \mathcal{Z}(s,\theta)$ has only one simple zero which is real. The proof is two-fold:
first we analyze zeros for bounded values of $\Im(s)$ (the low frequency part) using standard transfer operator techniques that go back to Parry-Pollicott 
\cite{ParryPollicott1,ParryPollicott2} and a Rouch\'e type argument of complex analysis. The high frequency part is harder and requires to recreate the arguments of \cite{Naud2}
and check that one can deal uniformly with the extra oscillating term coming from the Homology, as was done in \cite{OhWinter} for the congruence cocycle. We will choose a different (and faster) route by using a powerful new result of Bourgain-Dyatlov \cite{BD1} that allows to estimate certain oscillatory integrals with respect to Patterson-Sullivan measure.

\noindent {\bf Acknowledgements}. This work started while attending the workshop "emerging topics: quantum chaos and fractal uncertainty principle" at Princeton
IAS in october 2017, and the author thanks the organizers for this opportunity. FN is supported by ANR GeRaSic and IUF.

\section{Abelian covers and zeta factorization}
 \subsection{Structure of Abelian covers} 
 \begin{center}
\includegraphics[scale=0.45]{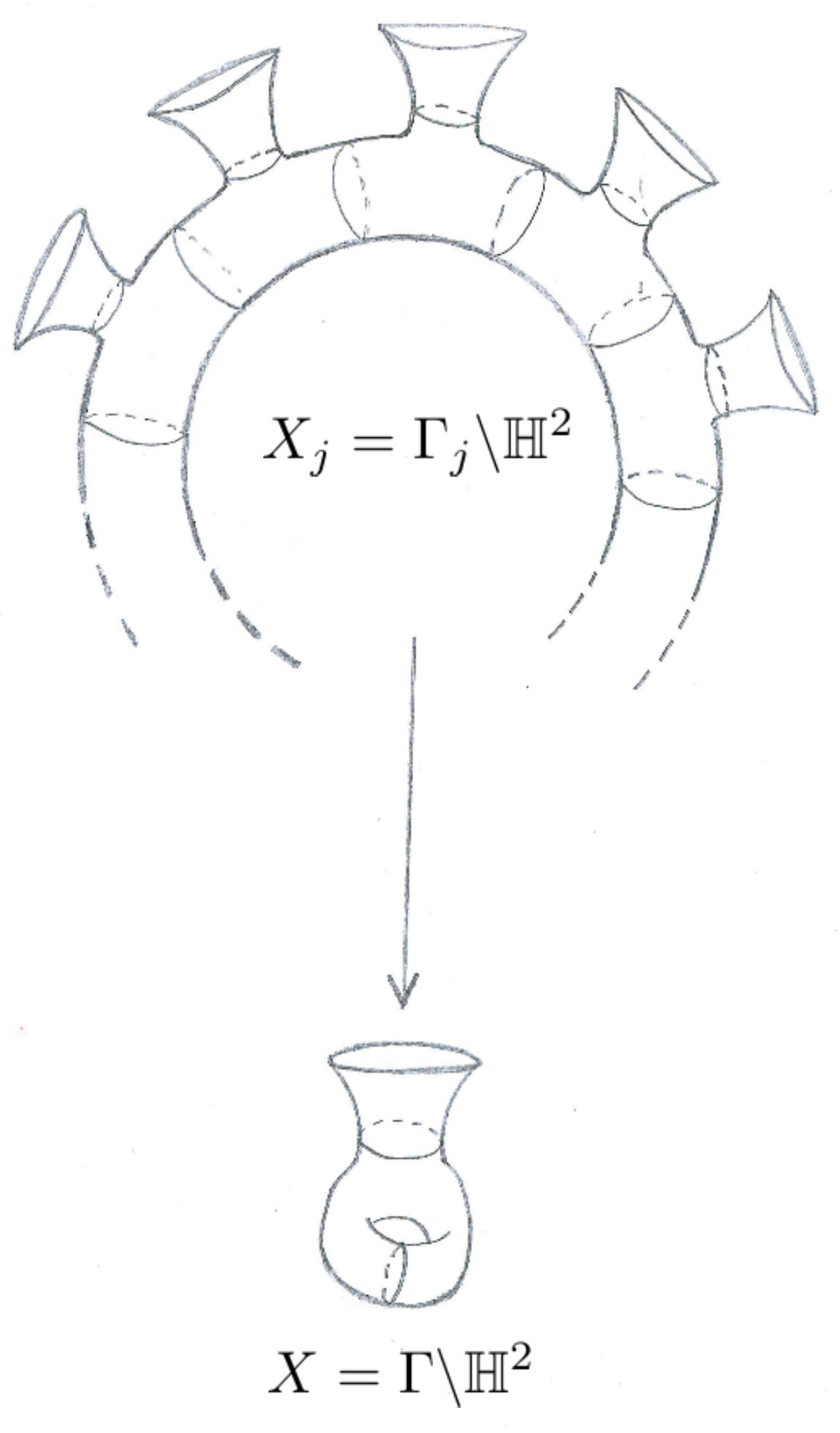} \\
\end{center} 
 Let $\Gamma$ be a convex co-compact group, then $\Gamma$ is isomorphic to the free group of rank $r$, with $r\geq 2$ when it is non elementary, see for example in \cite{Button} for a Schottky realization. Assume now that $\Gamma_j$ is a normal subgroup of $\Gamma$ such that
 $\G_j:=\Gamma/\Gamma_j$ is a finite abelian group. Let $\pi_j:\Gamma\rightarrow \G_j$ be the associated onto homomorphism so that 
 $$\Gamma_j=\ker(\pi_j).$$ By universal property of the abelianized group 
 $$\Gamma^{ab}:=\Gamma/[\Gamma,\Gamma]=H^1(X,\Z)\simeq \Z^r,$$
 the homomorphism $\pi_j$ can be factorized as $\pi_j=\widetilde{\pi_j}\circ P$ where $P:\Gamma\rightarrow \Z^r$ is a now fixed surjective homomorphism, 
 and $\widetilde{\pi_j}:\Z^r\rightarrow \G_j$ is another ($j$-dependent) onto homomorphism. By the usual structure theorem for finite Abelian groups, $\G_j$ can be written as
 a product of cyclic groups which we will write as 
 $$\G_j=\Z/N_1^{(j)}\Z\times \ldots\times  \Z/N_k^{(j)} \Z,$$
 where $N_1(j),N_2(j),\ldots N_k(j)$ are integers. 
 
 In the latter, we will be assuming that $1\leq k\leq r$ and that 
 $$\lim_{j\rightarrow +\infty} \inf_{\ell=1,\ldots,k} N_\ell^{(j)}=+\infty.$$
 Furthermore, $\widetilde{\pi_j}$ will be given by 
 $$\widetilde{\pi_j}(n)=(n_1\ \mathrm{mod}\ N_1^{(j)},\ldots,n_k\ \mathrm{mod}\ N_k^{(j)}),$$
 which is an obvious family of surjective homomorphism from $\Z^r$ to $\G_j$.

 In the simplest case $k=1$, the Galois group is the cyclic group $\Z/N_1^{(j)}\Z$, see the figure for an example, where
 the cover is obtained by cutting $X$ along a simple closed geodesic and glueing cyclically several copies of the result.
\subsection{Selberg's zeta function and characters}
According to the result of Patterson-Perry \cite{PatPerry}, resonances on $X=\Gamma \backslash \hh$ coincide with multiplicity with the non trivial zeros of the Selberg
zeta function, see also \cite{Borthwick} for the case of surfaces. Let $\mathcal{P}=\mathcal{P}(\Gamma)$ denote the set of primitive closed geodesics on $X$, and if $\mathcal{C}\in \mathcal{P}$,
$l(\mathcal{C})$ will be the length. Selberg zeta function is usually defined by the infinite product
$$Z_\Gamma(s):=\prod_{\mathcal{C}\in \mathcal{P}} \prod_{k \in \N_0}\left ( 1-e^{-(s+k)l(\mathcal{C})}\right),\ \Re(s)>\delta(\Gamma).$$
This infinite product has a holomorphic extension to $\C$. The characters of the abelian group 
$$H^1(X,\Z)\simeq \Z^r$$ 
are given by
$$\chi_\theta(x)=e^{2i\pi \langle \theta,x\rangle},\ x\in \Z,$$
where $\langle \theta,x\rangle=\sum_{\ell=1}^r \theta_\ell x_\ell$, and  $\theta=(\theta_1,\ldots,\theta_r)$ belongs to the torus
$\R^r/\Z^r$. Associated to each character $\chi_\theta$ is a corresponding "twisted" Selberg zeta $\mathcal{Z}_\Gamma(s,\theta)$ function (or rather $L$-function) defined by 
$$\mathcal{Z}_\Gamma(s,\theta):=\prod_{\mathcal{C}\in \mathcal{P}} \prod_{k \in \N_0}
\left ( 1-\chi_\theta(\mathcal{C})e^{-(s+k)l(\mathcal{C})}\right),\ \Re(s)>\delta(\Gamma),$$
where $\chi_\theta(\mathcal{C})$ is a shorthand for $\chi_\theta(P(\mathcal{C}))$. On the other hand, the characters of $\G_j$ are given by
$\chi_\theta((m_1,\ldots,m_k,0,\ldots,0)),\ m\in \G_j$, where
$$\theta \in \mathcal{S}_j:=\left \{0,\frac{1}{N_1^{(j)}}\ldots,\frac{N_1^{(j)}-1}{N_1^{(j)}}\right \} \times \ldots \times
\left \{0,\frac{1}{N_k^{(j)}}\ldots,\frac{N_k^{(j)}-1}{N_k^{(j)}}\right \}\times \underbrace{\{0\}\times \ldots \times \{0\}}_{r-k\  \mathrm{times} }.$$
Notice that if $\gamma \in \Gamma$, then for all $\theta \in \mathcal{S}_j$, we have indeed
$$\chi_\theta(\widetilde{\pi_j}\circ P(\gamma))=\chi_\theta(P(\gamma)).$$
From the results of \cite{JNS}, $\S 2$ and also \cite{PohlFedosova},
we know that for all $\theta \in \mathcal{S}_j$, each zeta function $s\mapsto \mathcal{Z}_\Gamma(s,\theta)$ has an analytic continuation to $\C$, and we have the following fundamental factorization formula, valid for all $s\in \C$:

\begin{equation}
\label{factor1}
Z_{\Gamma_j}(s)=\prod_{\theta \in\mathcal{S}_j} \mathcal{Z}_\Gamma(s,\theta).
\end{equation}
The main result then follows from the next Theorem.
\begin{thm}
\label{theo2}
Assume that $\Gamma$ is non-elementary. We have the following facts.
\begin{enumerate}
\item For all $\varepsilon>0$, one can find $\eta(\varepsilon)>0$ such that if $\theta \in \R^r$ is such that
$\mathrm{dist}(\theta,\Z^r)>\varepsilon$, then $s\mapsto \mathcal{Z}_\Gamma(s,\theta)$ does not vanish inside the strip
$$\{ \delta-\eta\leq \Re(s) \leq \delta\}.$$
\item There exists $\epsilon_0>0$ and $\eta_0>0$ such that for all $\theta$ with $\mathrm{dist}(\theta,\Z^r)\leq \epsilon_0$, the analytic function
$s\mapsto \mathcal{Z}_\Gamma(s,\theta)$ has exactly one zero $\varphi(\theta)$ (which is real) inside the strip
$$\{ \delta-\eta_0\leq \Re(s) \leq \delta\},$$
and the map $\theta\mapsto \varphi(\theta)$ is smooth, real valued with a non degenerate critical point at $\theta=0$. 
\end{enumerate}
\end{thm}
The proof of Theorem \ref{theo2} will occupy several sections. Let us show how one can recover Theorem \ref{theo1} from that. We first start by picking $\epsilon_0$ from 
statement (2), and then a corresponding $\eta(\epsilon_0)$ from statement (1). Set $\eta^*=\min\{ \eta_0;\eta(\epsilon_0)\}$. Inside the strip
$$\Omega:=\{ \delta-\eta^*\leq \Re(s) \leq \delta\}, $$
we observe that either $\mathrm{dist}(\theta,\Z^r)\leq \epsilon_0$ and $s\mapsto \mathcal{Z}_\Gamma(s,\theta)$ vanishes at most once on the real line,
or $\mathrm{dist}(\theta,\Z^r)> \epsilon_0$ and $s\mapsto \mathcal{Z}_\Gamma(s,\theta)$ does not vanish. Going back to the factorization formula (\ref{factor1}), we deduce that inside $\{ \delta-\eta^*\leq \Re(s) \leq \delta\},$ the set of zeros of $Z_{X_j}(s)$ is given by 
$$\{\varphi(\theta)\ :\ \theta \in \mathcal{S}_j\ \mathrm{and}\ \mathrm{dist}(\theta,\Z^r)\leq \epsilon_0 \} \cap  
\{ \delta-\eta^*\leq \Re(s) \leq \delta\}.$$
To complete the proof, we follow the arguments of \cite{JNS}, which we briefly recall for completeness. Let $f\in C_0^\infty([\delta-\epsilon_1,1])$,
where $0<\epsilon_1<\eta^*$ is small enough such that $\mathrm{Supp}(f\circ\varphi)\subset \{ \mathrm{dist}(\theta,\Z^r)\leq \epsilon_0\}$.
We therefore have
$$\frac{1}{\vert \G_j\vert} \sum_{\lambda \in \mathcal{R}_{X_j}\cap\Omega}f(\lambda)= 
\frac{1}{N_1^{(j)}\ldots N_k^{(j)}}\sum_{\beta \in \Z^k} f\circ\varphi \left(\frac{\beta_1}{N_1^{(j)}},\ldots,\frac{\beta_k}{N_k^{(j)}},0,\ldots,0\right).$$
Applying Poisson summation formula, we obtain that as $j\rightarrow +\infty$, the righthand side converges to
$$\int_{\R^k}f\circ \varphi(x,0,\ldots,0)dx=:\int f d\mu.$$
The fact that the push-forward measure $\mu$ is absolutely continuous follows from the non-degeneracy of the critical point at $0$, see \cite{JNS}. By further shrinking the strip (i.e. taking a smaller $\eta^*$), and a standard approximation argument, the proof 
of the first two claims is complete. We now prove the last point. First we observe that using Theorem 1.1 from \cite{JNS}, part (2), we have the existence of a constant $C_\Gamma>0$ such that for all $j$ and $s\in \C$, we have
\begin{equation}
\label{bound1}
\vert Z_{\Gamma_j}(s)\vert \leq C_\Gamma \exp\left( C_\Gamma \vert \G_j\vert \vert s\vert^2  \right).
\end{equation}
On the other hand, for all $\Re(s)>\delta$ and $\theta \in \mathcal{S}_j$, we have
$$\mathcal{Z}_\Gamma(s,\theta)=\exp\left( -\sum_{n=1}^\infty\frac{1}{n} \sum_{\mathcal{C}\in \mathcal{P}(X)}
\chi_\theta(\mathcal{C}^n)\frac{e^{-snl(\mathcal{C})}}{1-e^{-nl(\mathcal{C})} }\right),$$
which combined with the factorization formula (\ref{factor1}) shows that for $\Re(s)>\delta$, 
\begin{equation}
\label{bound2}
\vert Z_{\Gamma_j}(s)\vert\geq \exp\left( -C_1\vert \G_j\vert \sum_{n=1}^\infty\frac{1}{n} \sum_{\mathcal{C}\in \mathcal{P}(X)}
e^{-\Re(s)nl(\mathcal{C})}\right).
\end{equation}
We now fix $\lambda \in \mathcal{R}_X$ and $\varepsilon_0>0$. To get the upper bound we fix $x_0\in \R$ with 
$x_0>\delta$ and choose $R_0>0$ large enough such that the disc $D(x_0,R_0)$ contains $D(\lambda,\varepsilon_0)$ in its
interior. We will use Jensen's formula (or rather a consequence of it) in the following form. 
\begin{propo}
\label{Jensen}
 Let $f$ be a holomorphic function on the open disc $D(w,R)$, and assume that $f(w)\neq 0$. let $N_f(r)$ denote the number of zeros of $f$ in the closed disc $\overline{D}(w,r)$. For all $\widetilde{r}<r<R$, we have 
$$N_f(\widetilde{r})\leq \frac{1}{\log(r/\widetilde{r})} \left (  \frac{1}{2\pi} \int_0^{2\pi} \log \vert f(w+re^{i\theta})\vert d\theta-\log\vert f(w)\vert   \right).$$
\end{propo}
It is now clear that by applying the above Proposition on the disc $D(x_0,R_0)$ where both $x_0,R_0$ are fixed we can use
the bounds (\ref{bound1}), and (\ref{bound2}) to obtain that for all $j$,
$$\# \mathcal{R}_{X_j} \cap D(\lambda,\varepsilon_0)\leq C_\Gamma \vert \G_j \vert.$$
To prove the lower bound, provided $\varepsilon_0$ is taken small enough, we can write for all $s\in D(\lambda,\varepsilon_0)$,
$$Z_\Gamma(s)=(s-\lambda)^m\psi(s),$$
where $m\geq 1$ is the order of vanishing of $Z_\Gamma(s)$ at $s=\lambda$ and $s\mapsto \psi(s)$ is a holomorphic function non vanishing on a neighborhood of $\overline{D(\lambda,\varepsilon_0)}$. On $\partial D(\lambda,\varepsilon_0)$ we have
$$\vert Z_\Gamma(s)\vert \geq \epsilon_0^m \inf_{s\in D(\lambda,\varepsilon_0)} \vert \psi(s)\vert >0.$$
On the otherhand, since $(s,\theta)\mapsto \mathcal{Z}(s,\theta)$ is smooth and $\mathcal{Z}(s,0)=Z_\Gamma(s)$, there exist
$\epsilon>0$ such that for all $\Vert \theta \Vert\leq \epsilon$ we have
$$\sup_{s\in \partial D(\lambda,\varepsilon_0)}\vert \mathcal{Z}(s,\theta)-Z_\Gamma(s)\vert < 
\inf_{s\in \partial D(\lambda,\varepsilon_0)} \vert Z_\Gamma(s)\vert.$$
Applying the classical Rouch\'e's theorem for holomorphic functions, 
we deduce that for each $\theta\in \mathcal{S}_j$ such that $\Vert \theta \Vert \leq \epsilon$,
$s\mapsto \mathcal{Z}(s,\theta)$ has exactly $m$ zeros inside $D(\lambda,\varepsilon_0)$. Using the factorization formula,
we deduce that the number of zeros of $Z_{\Gamma_j}(s)$ inside $D(\lambda,\epsilon_0)$ is at least
$$m\#\{ \theta \in \mathcal{S}_j\ :\ \Vert \theta \Vert \leq \epsilon\},$$
which is bigger than $C\vert \G_j \vert$ for some small constant $C>0$, independent of $j$. The proof is complete.

\subsection{Closed geodesics in homology classes}
let $P:\Gamma \rightarrow \Z^r\simeq H^1(X,\Z)$ be a fixed isomorphism as above.
Let $\alpha \in \Z^r$ be a fixed "holomogy class", and consider the counting function
$$N(\alpha,T)=\#\{ \mathcal{C}\ \in \mathcal{P}(X)\ :\ P(\mathcal{C})=\alpha\ \mathrm{and}\ l(\mathcal{C})\leq T\}.$$
Counting asymptotics for closed geodesics in homology classes has a long history of results for compact hyperbolic manifolds or more general Anosov flows on compact manifolds, see \cite{Anantharaman,PhillipsSarnak,Pollicott2,Lalley2,KatSun}. In the case of infinite volume hyperbolic surfaces, the leading term is known, and follows for example from \cite{ParryPollicott2}, Chapter 12 (For Kleinian groups, we also mention the work of Babillot-Peign\'e \cite{BP1}). It goes as follows: as $T\rightarrow +\infty$ we have
\begin{equation}
\label{PNT}
N(\alpha,T)\sim c_0\frac{e^{\delta T}}{T^{r/2+1}},
\end{equation}
where $c_0$ is independent of $\alpha$.

As a consequence of Theorem \ref{theo2} on the non-vanishing of $\mathcal{Z}_\Gamma(s,\theta)$ and combining it with a priori estimates on zeta functions from \cite{JNS}, Theorem 1.1, we obtain the following improved counting result. 
\begin{thm}
Assume that $\Gamma$ is convex co-compact and non elementary, then for all $\alpha \in \Z^r$, for all $n\geq 0$, there exists a sequence
$c_0,c_1(\alpha),\ldots,c_n(\alpha)\in \R$ such that as $T\rightarrow +\infty$,
$$N(\alpha,T)=\frac{e^{\delta T}}{T^{r/2+1}}\left (c_0+c_1T^{-1}+\ldots+c_n T^{-n}+O(T^{-n-1}) \right).$$
\end{thm}
In particular, this extends the asymptotics obtained by McGowan and Perry \cite{McP} to the case $\delta\leq \half$, which was not known so far.
The proof, knowing Theorem \ref{theo2}, is standard and goes exactly as in \cite{McP}. We recall briefly the main ideas for the benefit of the reader. One starts by picking $\phi_T \in C_0^\infty(\R^+)$, $\phi_T\geq 0$ such that $\phi_T\equiv 1$ on the interval $[\epsilon_0,T]$ and is supported in $[\epsilon_0/2,T+\beta]$, where $\epsilon_0>0$ is taken small and $\beta=e^{-\nu T}$
for some large $\nu>0$. We then set
$$\psi_T(s):=\int_0^\infty e^{xs}\phi_T(x)dx,$$
so that for all $A>\delta$ we have the contour integral identity
$$\frac{1}{2i\pi} \int_{A-i\infty}^{A+i\infty} \frac{\mathcal{Z}'_\Gamma(s,\theta)}{\mathcal{Z}_\Gamma(s,\theta)}\psi_T(s)ds=
\sum_{k,\mathcal{C}} l(\mathcal{C})\frac{\chi_\theta(\mathcal{C}^k)}{1-e^{-kl(\mathcal{C})}}\phi_T(kl(\mathcal{C})).$$
Notice that if $\nu$ is large enough and $\epsilon_0$ small, we have for $\sigma\leq \delta$,
$$\phi_T(\sigma)=\frac{e^{\sigma T}}{\sigma}+O\left(e^{T\delta/2}\right). $$
Thanks to the a priori upper bound from \cite{JNS} and Caratheodory estimates, we know that if $\mathcal{Z}_\Gamma(s,\theta)\neq 0$ for all $s$ with $\Re(s)>\delta-\eta$,
then we will get a polynomial upper bound for the log derivative
$$ \left \vert \frac{\mathcal{Z}'_\Gamma(s,\theta)}{\mathcal{Z}_\Gamma(s,\theta)} \right \vert \leq M\vert \Im(s)\vert^2,$$
for all $\vert \Im(s)\vert$ large and $\Re(s)>\delta-\eta/2$.
Integrating with respect to $\theta$ on $\R^r/\Z^r$ gives the formula
$$\frac{1}{2i\pi} \int_{A-i\infty}^{A+i\infty} \int_{\R^r/\Z^r}e^{-2i\pi\langle \alpha,\theta\rangle} \frac{\mathcal{Z}'_\Gamma(s,\theta)}{\mathcal{Z}_\Gamma(s,\theta)}d\theta\psi_T(s)ds=
\sum_{k,\mathcal{C}\ :\ P(\mathcal{C}^k)=\alpha} \frac{l(\mathcal{C})}{1-e^{-kl(\mathcal{C})}}\phi_T(kl(\mathcal{C})).$$
Thanks to Theorem \ref{theo2}, for all $\epsilon>0$, we can therefore deform the contour (by taking $A<\delta$) for all $\theta$ such that $\mathrm{dist}(\theta,0)>\epsilon$
to obtain a contribution of order
$$O(e^{(\delta-\eta(\epsilon)/2)T}).$$
We are essentially left with estimating integrals over $\theta$ in a neighborhood of $0$. Using the {\it residue formula}, the fact that $s\mapsto \mathcal{Z}_\Gamma(s,\theta)$
has a simple leading zero $\varphi(\theta)$, and neglecting error terms which are exponentially smaller than $e^{\delta T}$, we are then led to estimate integrals of the form
$$\mathcal{I}(T)=\int_{\R^r/\Z^r} e^{\varphi(\theta)T}\kappa(\theta)d\theta,$$
where $\kappa(\theta)$ is a smooth function supported in an arbitrarily small neighborhood of $0$. Using Morse Lemma and Laplace method to deal with the {\it stationnary phase} at $\theta=0$ (see Lemma 2.3 in \cite{PhillipsSarnak})
leads to expansions as $T\rightarrow +\infty$ of the form
$$\mathcal{I}(T)=\frac{e^{\delta T}}{T^{r/2}}\left (a_0+a_1T^{-1}+\ldots+a_n T^{-n}+O(T^{-n-1}) \right). $$
Notice that there are no odd powers of $T^{-\half}$ here because all the odd moments on $\R^r$ of $e^{-\vert x\vert^2}$ vanish.
We have essentially obtained that
$$\sum_{P(\mathcal{C})=\alpha\ \mathrm{and}\ l(\mathcal{C})\leq T} l(\mathcal{C})=
\frac{e^{\delta T}}{T^{r/2}}\left (c_0+c_1T^{-1}+\ldots+c_n T^{-n}+O(T^{-n-1}) \right).$$
To obtain the desired asymptotics for $N(\alpha,T)$ is now a simple exercise using Stieltjes integration by parts and the bound coming from the known leading term (\ref{PNT}).
We point out that using more delicate arguments involving the saddle point method, it is possible to derive similar asymptotics for counting functions of the type
$$N(\alpha+[T\xi],T),$$
where $\xi\in \Z^r\setminus \{0\}$, see Anantharaman \cite{Anantharaman}.

\section{Schottky uniformization and transfer operators}

\subsection{Schottky groups}
We start by recalling Bowen-Series coding and holomorphic function spaces needed for our analysis. Let $\hh$ denote the Poincar\'e upper half-plane 
$$\hh=\{ x+iy\in \C\ :\ y>0\}$$
endowed with its standard metric of constant curvature $-1$
$$ds^2=\frac{dx^2+dy^2}{y^2}.$$ 
The group of isometries of $\hh$ is $\mathrm{PSL}_2(\R)$ through the action of 
$2\times 2$ matrices viewed as M\"obius transforms
$$z\mapsto \frac{az+b}{cz+d},\ ad-bc=1.$$ 
Below we recall the definition of Fuchsian Schottky groups which will be used to define transfer operators.
A Fuchsian Schottky group is a free subgroup of $\mathrm{PSL}_2(\R)$ built as follows. Let $\D_1,\ldots, \D_r,\D_{r+1},\ldots, \D_{2r}$, $r\geq 2$, 
be $2r$ Euclidean {\it open} discs in $\C$ orthogonal to the line $\R\simeq \partial \hh$. We assume that for all $i\neq j$, $\overline{\D_i} \cap \overline{\D_j}=\emptyset$. 
Let $\gamma_1,\ldots,\gamma_r \in \mathrm{PSL}_2(\R)$ be $r$ isometries such that for all $i=1,\ldots,r$, we have
$$\gamma_i(\D_i)=\widehat{\C}\setminus \overline{\D_{r+i}},$$
where $\widehat{\C}:=\C\cup \{ \infty \}$ stands for the Riemann sphere. For notational purposes, we also set $\gamma_i^{-1}=:\gamma_{r+i}$.
\begin{center}
\includegraphics[scale=0.65]{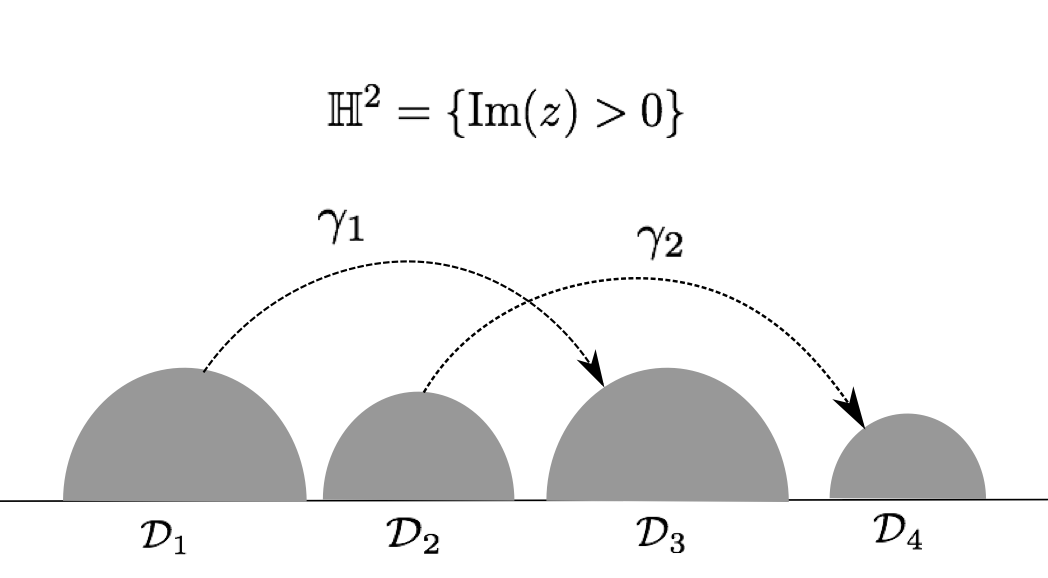}
\end{center}

\bigskip \noindent
Let $\Gamma$ be the free group generated by $\gamma_i,\gamma_i^{-1}$ for $i=1,\ldots,r$, then $\Gamma$ is a convex co-compact group, i.e. it is finitely generated and has no non-trivial parabolic element. The {\it converse is true} : up to isometry, convex co-compact hyperbolic surfaces
can be obtained as a quotient by a group as above, see \cite{Button}. 

For all $j=1,\ldots,2r$, set $I_j:=\D_j\cap \R$. One can define a map 
$$T:I:=\cup_{j=1}^{2r}I_j\rightarrow\R\cup\{\infty\}$$
by setting
$$T(x)=\gamma_j(x)\ \mathrm{if}\ x\in I_j.$$
This map encodes the dynamics of the full group $\Gamma$, and is called the Bowen-Series map, see 
\cite{BowenSeries1} for the genesis of these type of coding. The key properties are orbit equivalence and
uniform expansion of $T$ on the maximal invariant subset $\cap_{n\geq 1} T^{-n}(I)$ which coincides with the limit set
$\Lambda(\Gamma)$, see for example \cite{Borthwick}.

We now define the function space and the associated transfer operators.
Set 
$$\Omega:=\cup_{j=1}^{2r}\D_j.$$
Consider now the Hilbert space $H^2(\Omega)$ which is defined as the set of {\it holomorphic functions} 
$F:\Omega\rightarrow V_\rr$ such that
$$\Vert F\Vert_{H^2}^2:=\int_{\Omega} \vert F(z)\vert^2 dm(z)<+\infty, $$
where $dm$ is Lebesgue measure on $\C$. Let $\theta \in \R^r/\Z^r$, the "character torus". On the space $H^2(\Omega)$, we define a "twisted" by $\theta$ transfer operator $\lt_{s,\theta}$ by
$$\lt_{s,\theta}(F)(z):=
\sum_{j\neq i} (\gamma_j')^s(z) \chi_\theta(P\gamma_j)F(\gamma_j z),\ \mathrm{if}\ z\in \D_i,$$
where $s\in \C$ is the spectral parameter, and $\chi_\theta$ is the character of $H^1(X,\Z)\simeq \Z^r$ associated to $\theta$ and $P:\Gamma\rightarrow H^1(X,\Z)$
is the projection homomorphism.
Notice that for all $j\neq i$, $\gamma_j:\D_i\rightarrow \D_{r+j}$ is a holomorphic contraction since $\overline{\gamma_j(\D_i)}\subset \D_{r+j}$.
Therefore, $\lt_{s,\theta}$ is a compact {\it trace class} operator and thus has a {\it Fredholm determinant}. We define the twisted zeta function 
$\mathcal{Z}_\Gamma(s,\theta)$ by
$$\mathcal{Z}_\Gamma(s,\theta):=\det(I-\lt_{s,\theta}).$$
It follows from \cite{JNS}, but also \cite{PohlFedosova} that for all $\Re(s)>\delta$ we have the identity
$$ \det(I-\lt_{s,\theta})=\prod_{\mathcal{C}\in \mathcal{P}} \prod_{k \in \N_0}
\left ( 1-\chi_\theta(\mathcal{C})e^{-(s+k)l(\mathcal{C})}\right),$$
which shows that the infinite product has actually an analytic continuation to $\C$.

\subsection{The high-low frequency results}
The proof of Theorem \ref{theo2} will follow from two facts which will require two different types of asymptotic analysis. We state these results below.
\begin{propo}
\label{HighF}(The high frequency regime) Assume that $\Gamma$ is non elementary, then there exist $\varepsilon_0>0$ and $T_0>>1$ such that for
all $\theta \in \R^r$ and 
$$s\in \{\delta-\varepsilon_0\leq \Re(s)\leq \delta\ \mathrm{and}\ \vert \Im(s)\vert\geq T_0 \},$$
we have $\mathcal{Z}_\Gamma(s,\theta)\neq 0$.
\end{propo}
A very important feature is that $\varepsilon_0>0$ and $T_0$ can be taken {\it uniform} with respect to $\theta$. 
This uniform high frequency (aka large $\Im(s)$) fact will follow from
certain Dolgopyat estimates for twisted transfer operators as in \cite{Naud2}. 
In particular, this result implies that at high frequencies, there is a uniform resonance gap for all abelian covers
of a given non elementary Schottky surface, a fact that is similar to the result proved in \cite{OhWinter} for congruence subgroups.
To describe the behaviour of resonances with small $\Im(s)$, we will prove the following result.
\begin{propo}
 \label{LowF}(The low frequency regime) Assume that $\Gamma$ is non elementary, then for all $t\in \R$ and $\theta \in \R^r/\Z^r$ we have
 $$\mathcal{Z}_\Gamma(\delta+it,\theta)=0 \Longleftrightarrow (t,\theta)=(0,0),$$
 where $0$ in the second factor is understood mod $\Z^r$.
\end{propo}
In other words, on the vertical line $\{\Re(s)=\delta\}$, the zeta function $\mathcal{Z}_\Gamma(s,\theta)$ vanishes only at $s=\delta$ when $\theta\in \Z^r$.
The proof will follow from convexity arguments in the analysis of transfer operators, as in previous works of Parry and Pollicott \cite{ParryPollicott2}.

To conclude this section, let us show how the combination of Proposition \ref{HighF} and Proposition \ref{LowF} does imply Theorem \ref{theo2}. First we fix $\epsilon>0$.
We know from Proposition \ref{HighF} that no zeta function $\mathcal{Z}_\Gamma(s,\theta)$ will vanish for $\delta-\varepsilon_0\leq \Re(s)\leq \delta$ and $\vert \Im(s) \vert \leq T_0$
regardless of the value of $\theta$. Assume that for all $\eta>0$, there exists  $\theta \in \R^r$ with $\mathrm{dist}(\theta,\Z^r)>\epsilon$ and there exists $s\in \C$ with
$\delta-\eta\leq \Re(s)\leq \delta$ and $\vert \Im(s)\vert \leq T_0$ such that $\mathcal{Z}_\Gamma(s,\theta)=0$. Then by compactness one construct a converging sequence
$(s_\ell,\theta_\ell)$ such that 
$$s_\infty:=\lim_{\ell\rightarrow +\infty} s_\ell \in \delta+i[-T_0,+T_0]$$ 
and $\theta_\infty:=\lim_{\ell\rightarrow +\infty} \theta_\ell$ satisfies 
$\theta_\infty \not \in \Z^r$. By continuity, we have
$$\mathcal{Z}_\Gamma(s_\infty,\theta_\infty)=0$$
which clearly contradicts Proposition \ref{LowF}. Therefore one can find $\widetilde{\eta}(\varepsilon)>0$ such that if $\theta \in \R^r$ is such that
$\mathrm{dist}(\theta,\Z^r)>\varepsilon$, then $s\mapsto \mathcal{Z}_\Gamma(s,\theta)$ does not vanish inside the rectangle
$$\{ \delta-\widetilde{\eta}\leq \Re(s) \leq \delta\ \mathrm{and}\ \vert \Im(s)\vert \leq T_0\}.$$
By taking $\eta=\min\{\varepsilon_0 ,\widetilde{\eta}\}$ we have proved part (1) of Theorem \ref{theo2}.

Let us consider the family of rectangles
$$R_{T_0,\eta}:=[\delta-\eta,\delta+\eta]+i[-T_0,+T_0].$$
Because we have $\mathcal{Z}_\Gamma(s,0)=Z_\Gamma(s)$ and $(s,\theta)\mapsto \mathcal{Z}_\Gamma(s,\theta)$ is smooth, there exists a constant $C_{T_0,\eta}>0$ such that
for all $\theta \in \R^r$ with $\Vert \theta \Vert\leq \epsilon_0$ we have for all $s\in R_{T_0,\eta}$,
$$\vert \mathcal{Z}_\Gamma(s,\theta)-Z_\Gamma(s) \vert \leq C_{T_0,\eta}\epsilon_0.$$
On  the other hand, since on the line $\{ \Re(s)=\delta\}$, $Z_\Gamma(s)$ vanishes only at $s=\delta$, with a simple zero, one can find $\eta_0>0$ small enough such that for all $s \in R_{T_0,\eta_0}$ one can write
$$Z_\Gamma(s)=(s-\delta)\psi(s),$$
where $\psi(s)$ is holomorphic in a neighbourhood of $R_{T_0,\eta_0}$ and does not vanish on $R_{T_0,\eta_0}$. For all $s\in \partial R_{T_0,\eta_0}$, we have
$$\vert Z_\Gamma(s)\vert\geq  \eta_0\inf_{R_{T_0,\eta_0}}\vert \psi(s)\vert=:M_{T_0,\eta_0}.$$
By choosing $\epsilon_0>0$ small enough we can make sure that $M_{T_0,\eta_0}>\epsilon_0 C_{T,\eta_0}$ so that we can apply {\it Rouch\'e's theorem} to conclude that
$\mathcal{Z}_\Gamma(s,\theta)$ has exactly one simple zero in $R_{T_0,\eta_0}$. By combining it with Proposition \ref{HighF}, we now know that provided $\Vert \theta \Vert$ is small enough, $s\mapsto \mathcal{Z}_\Gamma(s,\theta)$ has exactly one zero in a thin strip $\{ \delta-\eta_0\leq \Re(s)\leq \delta \}$. The fact that this zero is real follows from "time reversal"
invariance of the length spectrum: in other words, we have
$$\overline{\mathcal{Z}_\Gamma(s,\theta)}=\mathcal{Z}_\Gamma(\overline{s},\theta),$$
see \cite{JNS} for more details. Since non real zeros must come in conjugate pairs, this forces this unique zero to be real. The fact that this unique zero can be smoothly parametrized
as a function $\varphi(\theta)$ for all $\Vert \theta \Vert$ small is just an application of the holomorphic implicit function theorem, legitimate since
$$\partial_s \mathcal{Z}_\Gamma(\delta,0)=Z_\Gamma'(\delta)\neq 0.$$
Because we have $\varphi(0)=\delta$, and $\varphi(\theta)\leq \delta$ for all $\theta$ close to $0$ (indeed, all zeta functions $\mathcal{Z}_\Gamma(s,\theta)$ 
do not vanish inside $\{\Re(s)>\delta\}$),  the map $\theta\mapsto \varphi(\theta)$ must have a critical point at $s=\delta$. The fact that it is non degenerate is more subtle,
see \cite{JNS} for references.

\section{High frequency analysis and uniform Dolgopyat estimates}

The goal of this section is to prove Proposition \ref{HighF} which is concerned with zeros of $\mathcal{Z}_\Gamma(s,\theta)$ for $\Re(s)$ 
close to $\delta$ and large $\vert \Im(s)\vert$. When $\theta=0\ \mathrm{mod}\ \Z^r$, then this was done in \cite{Naud2}. The game here is to show
that one can do the same {\it uniformly} in $\theta$. As pointed out in \cite{OhWinter}, the fact that the extra character term $\chi_\theta(\gamma)$
is {\it locally constant} on $I=\cup_j I_j$ makes it possible to apply almost verbatim the analysis of \cite{Naud2}, where one has essentially to check that the extra oscillating term does not interfere with the "large $\Im(s)$" cancellation mechanism.

 In this section we will choose an {\it alternative route} based on the recent result of \cite{BD1} which will allow us to bypass the {\it most technical part} of the argument in \cite{Naud2}, allowing an easier proof of the {\it uniform spectral gap}. We believe
this alternative proof might be interesting for future generalizations of \cite{OhWinter} to arbitrary families of non-Galois covers, this will be pursued elsewhere.
\subsection{Reduction to $L^2$ estimates}
Let $C^1(I)$ denote the Banach space of complex valued functions, $C^1$ on $\overline{I}$, endowed with the norm ($t\neq 0$)
$$\Vert f\Vert_{(t)}:=\Vert f \Vert_\infty+\frac{1}{\vert t \vert} \Vert f' \Vert_\infty,$$
where as usual
$$\Vert f \Vert_\infty=\sup_{x \in I} \vert f(x) \vert.$$
We recall that the action of the transfer operator $\lt_{s,\theta}$, now on $C^1(I)$, is given by
$$\lt_{s,\theta}(F)(x):=
\sum_{j\neq i} (\gamma_j')^s(x) \chi_\theta(P\gamma_j)F(\gamma_j x),\ \mathrm{if}\ x\in I_i.$$
We need to recall a few basic estimates that we will use throughout the rest of the paper. 

We first introduce some notations. We recall that $\gamma_1,\ldots,\gamma_r$ are generators of the Schottky group $\Gamma$, as defined in the previous section.
Considering a finite sequence $\alpha$ with
\[\alpha=(\alpha_1,\ldots,\alpha_n)\in \{1,\ldots, 2r\}^n,\]
we set 
\[ \gamma_\alpha:=\gamma_{\alpha_1}\circ \ldots \circ \gamma_{\alpha_n}. \]
We then denote by $\mathscr{W}_n$ the set of admissible sequences of length $n$ by
\[ \mathscr{W}_n:=\{ \alpha \in \{1,\ldots, 2r\}^n\ :\ 
\forall\ i=1,\ldots,n-1,\ \alpha_{i+1}\neq \alpha_i +r\ \mathrm{mod}\ 2r \}.\]
We point out that $\alpha \in \mathscr{W}_n$ if and only if $\gamma_\alpha$ is a {\it reduced word}
in the free group $\Gamma$. For all $j=1,\ldots, 2r$, we define $\mathscr{W}_n^j$ by
\[ \mathscr{W}_n^j:=\{ \alpha \in \mathscr{W}_n\ :\ \alpha_n\neq j \}. \] 
If $\alpha \in \mathscr{W}_n^j$, then $\gamma_\alpha$ maps $\overline{\D_j}$ into $\D_{\alpha_1+r}$.
Given the above notations and $f\in C^1(I)$, we have for all $x\in I_j$ and $n\in \N$,
$$\lt_{s,\theta}^n(f)(x)=\sum_{\alpha \in \mathscr{W}_n^j} (\gamma_\alpha'(x))^s \chi_\theta(P\gamma_\alpha)f(\gamma_\alpha (x)).$$

We will need throughout the paper some distortion estimates
for these "inverse branches" of $T^n$ that can be found in \cite{Naud2}. More precisely we have:
\begin{itemize}
 \item {\it(Uniform hyperbolicity)}. One can find $C>0$ and $0<\overline{\rho}<\rho<1$ such that for all $n$ and all $j$ such that $\alpha \in \mathscr{W}_n^j$, then for all $x\in I_j$ we have
 \[ C^{-1}\overline{\rho}^n\leq \vert \gamma'_\alpha(x) \vert \leq C \rho^n. \]
 \item {\it (Bounded distortion).} There exists $M_1>0$ such that for all $n,j$ and all $\alpha \in \mathscr{W}_n^j$,
 \[ \sup_{I_j} \left \vert \frac{\gamma''_\alpha}{\gamma'_\alpha}  \right \vert \leq M_1.\]
\item {\it (Bounded distortion for third derivatives)}. There exists $Q>0$ such that for all $n,j$ and all $\alpha \in \mathscr{W}_n^j$,
 \[ \sup_{I_j} \left \vert \frac{\gamma'''_\alpha}{\gamma'_\alpha}  \right \vert \leq Q.\]
\end{itemize}
The bounded distortion estimate has the following important consequence: there exists a uniform constant $M_2>0$ such that
for all $x,y \in I_j$, 
$$\left \vert \frac{\gamma'_\alpha(x)}{\gamma'_\alpha(y)}  \right \vert \leq M_2. $$
We point out that the same conclusion is still valid (up to a bigger constant $M_3$) if $x$ and $y$ belong to different $I_j$ and
$I_{j'}$ such that $\alpha \in  \mathscr{W}_n^j\cap \mathscr{W}_n^{j'} $. Indeed if $\alpha=\alpha_1\ldots \alpha_n$ then both
$\gamma_{\alpha_n}(x),\gamma_{\alpha_n}(y)\in I_{\alpha_n+r}$ and we can apply the above estimate.

\noindent
We will also need to use the {\it topological pressure} and {\it Bowen's formula}. 
Recall that the Bowen-Series map  $$T:\cup_{i=1}^{2p} I_i\rightarrow \R\cup\{\infty \}$$ is defined by $T(x)=\gamma_i(x)$ if $x\in I_i$. The non-wandering set of this map is exactly the limit set 
$\Lambda(\Gamma)$ of the group:
\[ \Lambda(\Gamma)=\bigcap_{n=1}^{+\infty} T^{-n}(\cup_{i=1}^{2p} I_i).\]
The limit set is $T$-invariant and given a continuous map $\varphi:\Lambda(\Gamma)\rightarrow \R$, the topological pressure $P(\varphi)$ can be defined through the variational formula:
\[ P(\varphi)=\sup_{\mu}\left ( h_\mu(T)+\int_{\Lambda} \varphi d\mu \right),\]
where the supremum is taken over all $T$-invariant probability measures on $\Lambda$, and $h_\mu(T)$ stands for the measure-theoretic entropy.
A celebrated result of Bowen \cite{Bowen1} says that the map 
$$\sigma\mapsto P(-\sigma \log  T' )$$
is convex, strictly decreasing and vanishes
exactly at $\sigma=\delta(\Gamma)$, the Hausdorff dimension of the limit set. An alternative way to compute the topological pressure is to look at weighted sums
on periodic orbits i.e. we have 
\begin{equation}
\label{pressure0}
e^{P(\varphi)}=\lim_{n\rightarrow +\infty} \left( \sum_{T^n x=x} e^{\varphi^{(n)}(x)}\right)^{1/n},
\end{equation}
with the notation $\varphi^{(n)}(x)=\varphi(x)+\varphi(Tx)+\ldots +\varphi(T^{n-1}x).$
We will use the following fact.
\begin{lem} 
\label{pressure1}
For all $\sigma_0,M$ in $\R$ with $0\leq \sigma_0<M$, one can find $C_0>0$ such that for all $n$ large enough and $M\geq \sigma\geq \sigma_0$, we have
\begin{equation}
\label{pressure2}
 \sum_{j=1}^{2p}\left (\sum_{\alpha \in \mathscr{W}_n^j} \sup_{I_j}   (\gamma'_\alpha  )^\sigma \right)\leq C_0 e^{nP(\sigma_0)},
\end{equation}
where $P(\sigma)$ is used as a shorthand for $P(-\sigma T')$.
\end{lem}
The proof of this Lemma follows rather straightforwardly from the Ruelle-Perron-Frobenius Theorem, which we state below (\cite{ParryPollicott1}, Theorem 2.2), and will be used several times.
\begin{propo}[Ruelle-Perron-Frobenius] Set $\lt_\sigma=\lt_{\sigma,0}$ where $\sigma$ is real. 
\begin{itemize}
\item The spectral radius of $\lt_\sigma$ on $C^1(I)$ is $e^{P(\sigma)}$ which is a simple eigenvalue associated to a strictly positive
eigenfunction $h_\sigma>0$ in $C^1(I)$.
\item The operator $\lt_\sigma$ on $C^1(I)$ is quasi-compact with essential spectral radius smaller than $\kappa(\sigma)e^{P(\sigma)}$ for some 
$\kappa(\sigma)<1$.
\item There are no other eigenvalues on $\vert z\vert=e^{P(\sigma)}$.
\item Moreover, the spectral projector $\mathbb{P}_\sigma$ on $\{e^{P(\sigma)}\}$ is given by
$$\mathbb{P}_\sigma(f)=h_\sigma \int_{\Lambda(\Gamma)} f d\mu_\sigma,$$
where $\mu_\sigma$ is the unique $T$-invariant probability measure on $\Lambda$ that satisfies
$$\lt_\sigma^*(\mu_\sigma)=e^{P(\sigma)}\mu_\sigma.$$
\end{itemize}
\end{propo}
We continue with a basic a priori estimate.
\begin{lem}[Lasota-Yorke estimate]
Fix some $\sigma_0<\delta$, then there exists $C_0>0,\ \rho<1$ such that for all $n,\theta$ and all $s=\sigma+it$ with $\sigma\geq \sigma_0$, we have
$$\Vert \left (\lt_{s,\theta}^n(f)\right)'\Vert_\infty\leq C_0e^{nP(\sigma_0)}\left \{ (1+\vert t \vert )\Vert f\Vert_\infty+\rho^n\Vert f' \Vert_\infty    \right\}.$$
\end{lem}
\noindent {\it Proof}. Differentiate the formula for $\lt^n_{s,\theta}(f)$ and then use the bounded distortion property plus the uniform contraction, combined with the pressure estimate (\ref{pressure2}). Uniformity with respect to $\theta$ follows from the fact that $\vert \chi_\theta\vert\equiv 1$. $\square$

The main result of this section is the following.
\begin{propo}[Uniform Dolgopyat estimate]
\label{Decay1} There exist $\epsilon>0$, $T_0>0$ and $C,\beta>0$ such that for all $\theta$ and $n=[C\log\vert t\vert]$ with $s=\sigma+it$ satisfying 
$\vert \sigma-\delta\vert\leq \epsilon$ and $\vert t\vert\geq T_0$, we have
$$\int_{\Lambda(\Gamma)} \vert \lt_{s,\theta}^n(f)\vert^2d\mu_\delta \leq \frac{\Vert f\Vert^2_{(t)}}{\vert t\vert^\beta}.$$
\end{propo}
This type of estimate is very similar in spirit to the ones encountered in the seminal work of Dolgopyat \cite{Dolg} on Anosov flows, hence the terminology.
We claim that Proposition \ref{Decay1} implies Proposition \ref{HighF}. Assume that $\sigma\leq \delta$. First we observe that if $g\in C^1(I)$ is {\it positive}, then we write ($x\in I_j$)
$$\lt_\sigma^n(g)(x)=\sum_{\alpha \in \mathscr{W}_n^j} (\gamma_\alpha'(x))^\sigma g(\gamma_\alpha (x)),$$
and using the uniform hyperbolicity estimate (the lower bound) we have
$$\lt_\sigma^n(g)(x)\leq A(\sigma,n)\lt^n_\delta(g),$$
where $A(\sigma,n)\leq C\overline{\rho}^{(\sigma-\delta)n}$. Now write $n=n_1+n_2$ where both $n_1,n_2$ will be specified later on. Given $f\in C^1(I)$, we write
$$\Vert \lt_{s,\theta}^n(f)\Vert_\infty\leq A(\sigma,n_1)\Vert \lt_\delta^{n_1}(\vert \lt_{s,\theta}^{n_2}(f)\vert)\Vert_\infty.$$
Using the Ruelle-Perron-Forbenius theorem at $\sigma=\delta$ gives (using Cauchy-Schwarz and the fact that $\mu_\delta$ is a probability measure)
$$\Vert \lt_{s,\theta}^n(f)\Vert_\infty\leq CA(\sigma,n_1)\left ( \left (\int_{\Lambda(\Gamma)} \vert \lt_{s,\theta}^{n_2}(f)\vert^2d\mu_\delta\right)^{1/2}+\kappa^{n_1}
\Vert \lt_{s,\theta}^{n_2}(f) \Vert_{C^1} \right),$$
for some $0<\kappa<1$. Using the Lasota-Yorke estimate, we know that for $\sigma_0\leq \sigma\leq \delta$ we have (assume $\vert t\vert \geq 1$)
$$\Vert \lt_{s,\theta}^{n_2}(f) \Vert_{C^1}\leq C_0e^{n_2P(\sigma_0)}\vert t \vert \Vert f\Vert _{(t)}.$$
Using Proposition \ref{Decay1} with $n_2=[C_2\log\vert t\vert ]$, we get for $\vert t\vert \geq T_0$ and $\sigma_0\leq \sigma\leq \delta$ with $\vert \sigma_0-\delta\vert \leq \epsilon$, 
$$\Vert \lt_{s,\theta}^n(f)\Vert_\infty\leq CA(\sigma_0,n_1)\left ( \frac{1}{\vert t\vert^{\beta/2}}+\kappa^{n_1}
\vert t\vert^{1+C_2P(\sigma_0)}  \right)\Vert f\Vert_{(t)}.$$
We know choose $n_1=[C_1\log\vert t\vert]$ with $C_1$ large enough so that for $\vert t\vert \geq T_0$, we have
$$\Vert \lt_{s,\theta}^n(f)\Vert_\infty\leq CA(\sigma_0,n_1)\frac{\Vert f\Vert_{(t)}}{\vert t\vert^{\beta/2}},$$
and since we have
$$A(\sigma_0,n_1)\leq C \vert t\vert^{\widetilde{C_1}(\delta-\sigma_0)},$$
with $\widetilde{C_1}=C_1\vert\log \overline{\rho} \vert$, we can make sure that $\sigma_0$ is taken close enough to $\delta$ so that 
$$\Vert \lt_{s,\theta}^{n}(f)\Vert_\infty\leq \frac{\Vert f\Vert_{(t)}}{\vert t\vert^{\beta/4}},$$
for all $\vert t\vert\geq T_0$. Using the Lasota-Yorke estimate, a similar computation leads to the conclusion that for all $\theta$ and 
$n(t)=[C_3\log\vert t\vert]$ for some $C_3>0$, we get 
\begin{equation}
\label{Dolg1}
\Vert \lt_{s,\theta}^{n(t)}(f)\Vert_{(t)}\leq  \frac{\Vert f\Vert_{(t)}}{\vert t\vert^{\overline{\beta}}},
\end{equation}
for some $\overline{\beta}>0$ and $\vert t\vert \geq T_0>>1$, $\vert \sigma-\sigma_0\vert \leq \epsilon$ with $\epsilon>0$. 
Assume now that $\mathcal{Z}_\Gamma(s,\theta)$ has a zero inside the region 
$$\{s\in \C\ :\ \vert \Re(s)-\delta\vert\leq \epsilon\ \mathrm{and}\ \vert \Im(s)\vert \geq T_0\}.$$
Then we get the existence of $f_{s,\theta}\in C^1(I)$ with $\Vert f_{s,\theta}\Vert_{\vert\Im(s)\vert}=1$ such that
$$\lt_{s,\theta}(f_{s,\theta})=f_{s,\theta}.$$
Using (\ref{Dolg1}) this leads to
$$1 \leq \frac{1}{\vert \Im(s) \vert^{\overline{\beta}}},$$
which is clearly a contradiction since $\vert \Im(s) \vert >>1$.

The remaining subsections will focus on proving Proposition \ref{Decay1}.

\subsection{The measure $\mu_\delta$ versus Patterson-Sullivan density at $i$.}
Patterson-Sullivan densities are measures on the limit set that satisfy interesting invariance properties. In the convex co-compact case, they where first introduced by Patterson in \cite{Patterson}. Primarily defined on the disc model ${\mathbb D}$ of the hyperbolic plane, they are constructed via Poincar\'e series
(with $s>\delta(\Gamma)$, $x\in \mathbb{D}$)
$$P_\Gamma(s;x,x):= \sum_{\gamma \in \Gamma}e^{-s\mathrm{d}(x,\gamma x)}.$$
By taking weak limits as $s\rightarrow \delta$ of probability measures 
$$\nu_{x, s}:=\frac{\sum_{\gamma \in \Gamma}e^{-s\mathrm{d}(x,\gamma x)}\mathcal{D}_{\gamma x}}{P_\Gamma(s;x,x)},$$
where $\mathcal{D}_z$ is the Dirac mass at $z\in \mathbb{D}$, one obtains a $\Gamma$-invariant measure $\nu_x$ supported on the limit set. For the
upper half-plane model $\hh$, one can use the push forward of $\nu_x$ by the inverse of the Cayley map given by $A(z)=i\left (\frac{1+z}{1-z} \right)$.

The Patterson-Sullivan density $\nu:=\nu_i$ (centered at $i$) is then a probability measure supported on the limit set $\Lambda(\Gamma)\subset \R$
that satisfies the equivariant formula (for any integrable  $f$ on $\Lambda(\Gamma)$)
$$\forall\ \gamma \in \Gamma,\ \int_{\Lambda(\Gamma)}fd\nu= 
\int_{\Lambda(\Gamma)} (f\circ \gamma) \vert \gamma' \vert_{\mathbb D}^\delta d\nu,$$
Where $\vert \gamma'(x) \vert_{\mathbb D}$ comes from the unit disc model of $\hh$, given explicitly by
$$\vert \gamma'(x) \vert_{\mathbb D} :=\gamma'(x)\left (\frac{1+x^2}{1+\gamma(x)^2}\right ).$$
See for example Borthwick  \cite{Borthwick}, Lemma 14.2.
This Patterson-Sullivan density $\nu$ is actually absolutely continuous with respect to $\mu_\delta$, more precisely we have
the following.
\begin{lem}
\label{PS1}
There exists $C_\Gamma>0$ such that the measure $\mu_\delta$ from the Ruelle-Perron-Frobenius theorem is 
$$\mu_\delta=C_\Gamma{(1+x^2)^\delta}\nu.$$
\end{lem}
\noindent {\it Proof.} From the equivariant formula, we know that for all integrable $f$ and all  bounded interval $J$ we have
for all $\gamma \in \Gamma$,
$$\int_J f d\nu=\int_{\gamma^{-1}(J)} (f\circ \gamma) \vert \gamma' \vert_{\mathbb D}^\delta d\nu.$$
Remark that 
$$\Lambda(\Gamma)\subset \bigcup_{j=1,\ldots,2r}\bigcup_{i\neq j} \gamma_i(I_j),$$
so that we write
$$\int_{\Lambda(\Gamma)} f d\nu=
\sum_j \sum_{i\neq j} \int_{\gamma_i(I_j)} fd\nu.$$
By using the equivariant formula as above we get
$$\int_{\Lambda(\Gamma)} f d\nu=\sum_j  \int_{I_j} \sum_{i\neq j}(f\circ \gamma_i ) \vert \gamma'_i \vert_{\mathbb D}^\delta d\nu,$$
which we recognize as
$$\int_{\Lambda(\Gamma)} f d\nu=\int_{\Lambda(\Gamma)} H^{-1}(x)\lt_\delta(Hf)(x)d\nu(x),$$
where $H(x)=\frac{1}{(1+x^2)^\delta}$. It is now clear that acting on measures, we have
$$\lt^*_\delta( H^{-1}\nu)=H^{-1}\nu,$$
which by uniqueness of $\mu_\delta$ (normalized as a probability measure) implies the statement. $\square$

Since the density $H^{-1}$ is smooth and uniformly bounded from above and below on $\Lambda(\Gamma)$, the measure
$\mu_\delta$ inherits straightforwardly  some of the properties of Patterson-Sullivan densities. In particular we will need to use the following bound.
\begin{propo}[Ahlfors-David upper regularity] 
\label{ADR}
There exists $B_\Gamma>0$ such that for all bounded interval $J$, 
$$\mu_\delta(J)\leq B_\Gamma \vert J\vert^\delta.$$
\end{propo}
For a proof (for $\nu$) of that fact see for example \cite{Borthwick}, Lemma 14.13.
In \cite{BD1}, Bourgain-Dyatlov established the following {\it remarkable} Theorem.
\begin{thm}[Decay of oscillatory integrals]
\label{Decay2}
There exist constants  $\beta_1,\beta_2 >0$ such
that the following holds. Given $g\in C^1(I)$ and $\Phi \in C^2(I)$,  
consider the integral 
$$\mathcal{I}(\xi):=\int_{\Lambda(\Gamma)} e^{-i\xi \Phi(x)}g(x)d\nu(x).$$
If we have
$$\epsilon:=\inf_{\Lambda(\Gamma)} \vert \Phi'\vert >0,$$
and $\Vert \Phi \Vert_{C^2}\leq M$,
then for all $\vert \xi\vert \geq 1$, we have
$$\vert \mathcal{I}(\xi)\vert\leq C_{M} \vert \xi \vert^{-\beta_1} \epsilon^{-\beta_2}\Vert g \Vert_{C^1},$$
where $C_{M}>0$ does not depend on $\xi,\epsilon,g$.
\end{thm}
\noindent {\bf Remarks}. This result is stated as Theorem 2 in \cite{BD1}. However the dependence on $g$ and $\epsilon$ is not explicit in their statement.
The fact that it can be bounded using $\Vert g \Vert_{C^1}$ is obvious: it follows from linearity in $g$ and Banach-Steinhaus theorem. The dependence on $\epsilon$ appears only in Lemma 3.5 of \cite{BD1}, where one can check that the loss is polynomial in $\epsilon^{-1}$. All we need is to allow $\epsilon\geq \vert \xi \vert^{-\kappa}$ for some $\kappa>0$ without ruining
the decay in $\vert \xi \vert$, see $\S 4.4$. We mention the recent related work of Jialun Li \cite{Li}, where similar bounds are proved.

By Lemma \ref{PS1}, it is clear that the exact same statement holds for $\mu_\delta$. The proof of Proposition \ref{Decay1}
will follow rather directly from this decay result and some additional facts that we will prove below.
\subsection{A uniform non integrability (UNI) result}
Given two words $\alpha,\beta \in \mathscr{W}_n^j$, consider the quantity
$$\mathscr{D}(\alpha,\beta):=
\inf_{x\in I_j} \left \vert \frac{\gamma''_\alpha(x)}{\gamma'_\alpha(x)} -\frac{\gamma''_\beta(x)}{\gamma'_\beta(x)}\right \vert.$$
We prove the following estimate, which will be used when estimating the "near-diagonal" sums (see the next section below). This type of estimate is
a generalization to Schottky groups of the work done by Baladi and Vall\'ee for the Gauss map \cite{BaladiVallee}.
\begin{propo}[UNI]
\label{UNI1}
There exist constants $M>0$ and $\eta_0>0$ such that for all $n$ and all $\epsilon=e^{-\eta n}$ with $0<\eta<\eta_0$, we have
for all $\alpha \in \mathscr{W}_n^j$,
$$\sum_{\beta \in \mathscr{W}_n^j,\ \mathscr{D}(\alpha,\beta)<\epsilon} \Vert \gamma'_\beta \Vert_{I_j,\infty}^\delta\leq M \epsilon^\delta.$$
\end{propo}
\noindent {\it Proof}. First we set some notations. If $\alpha$ is an admissible word in say $\mathscr{W}_n^j$, we will write 
$$\gamma_\alpha(x)=\frac{a_\alpha x+b_\alpha}{c_\alpha x+d_\alpha},\ a_\alpha d_\alpha-b_\alpha c_\alpha=1.$$
Each $\gamma_\alpha$ is a hyperbolic isometry of $\hh$ whose attracting fixed point will be denoted by $x_\alpha$ and repelling by
$x_\alpha^*$. The isometric circle of $\gamma_\alpha$ is the circle centered at 
$$z_\alpha=-\frac{d_\alpha}{c_\alpha}=\gamma_\alpha^{-1}(\infty),$$ 
with radius $\frac{1}{c_\alpha}$. We point out that by our definition of Schottky groups, we must have 
\begin{equation}
\label{ineq1}
\vert \gamma_\alpha^{-1}(\infty)\vert \leq M,
\end{equation}
for some uniform $M>0$. Since $x_\alpha^*$ is in the disc centered at $-\frac{d_\alpha}{c_\alpha}$ and of radius $1/\vert c_\alpha \vert$, we have obviously
$$\left\vert x_\alpha^*+\frac{d_\alpha}{c_\alpha}\right \vert \leq \frac{1}{\vert c_\alpha \vert}.$$
On the other hand, since we have
$$\Im(\gamma_\alpha(i))=\frac{1}{c_\alpha^2+d_\alpha^2},$$
we can use (\ref{ineq1}) to deduce that
$$\frac{1}{\vert c_\alpha \vert}\leq \widetilde{M}\sqrt{\Im(\gamma_\alpha(i))}.$$
By the hyperbolicity estimate, it is now easy to see that one can find constants $M',\eta_0>0$ such that for all $n$ we have
$$\frac{1}{\vert c_\alpha \vert}\leq M'e^{-\eta_0 n },$$
which in turn implies
\begin{equation}
\label{ineq2}
\left\vert x_\alpha^*+\frac{d_\alpha}{c_\alpha}\right \vert \leq M'e^{-\eta_0 n }.
\end{equation}
This estimate just says that repelling fixed point and center of isometric circle are exponentially close when the word length goes to infinity, a quantitative version of the well known fact that centers of isometric circles accumulate on the limit set.

Given $\gamma_\alpha$, then $\gamma_\alpha^{-1}=\gamma_{\overline{\alpha}}$, where 
$\overline{\alpha}=(\alpha_n+r)\ldots(\alpha_1+r)$, understood mod $2r$. We will use below the fact that $x_\alpha^*=x_{\overline{\alpha}}$,
and that $\gamma'_\alpha(x_\alpha)=\gamma'_{\overline{\alpha}}(x_\alpha^*)$.
We now go back to the quantity
$$\sum_{\beta \in \mathscr{W}_n^j,\ \mathscr{D}(\alpha,\beta)<\epsilon} \Vert \gamma'_\beta \Vert_{I_j,\infty}^\delta. $$
For each $\beta$ as above,  write
$$ \frac{\Vert \gamma'_\beta \Vert_{I_j,\infty}}{\gamma'_{\overline{\beta}}(x_{\overline{\beta}})}=
\frac{\Vert \gamma'_\beta \Vert_{I_j,\infty}}{\gamma'_{\beta}(x_{\beta})}\frac{\gamma'_{\overline{\beta}}(x_{\overline{\beta}})}{\gamma'_{\overline{\beta}}(x_{\overline{\beta}})},$$
which by the bounded distortion estimate gives 
$$\sum_{\beta \in \mathscr{W}_n^j,\ \mathscr{D}(\alpha,\beta)<\epsilon} \Vert \gamma'_\beta \Vert_{I_j,\infty}^\delta
\leq M''\sum_{\beta \in \mathscr{W}_n^j,\ \mathscr{D}(\alpha,\beta)<\epsilon} ( \gamma'_{\overline{\beta}}(x_{\overline{\beta}}))^\delta.$$
Using the Gibbs property for the $\mu_\delta$ measure, see \cite{ParryPollicott2} Corollary 3.2.1, we obtain that 
$$ \sum_{\beta \in \mathscr{W}_n^j,\ \mathscr{D}(\alpha,\beta)<\epsilon} \Vert \gamma'_\beta \Vert_{I_j,\infty}^\delta
\leq C' \sum_{\beta \in \mathscr{W}_n^j,\ \mathscr{D}(\alpha,\beta)<\epsilon} \mu_\delta(I_{\overline{\beta}}),$$
where $I_{\overline{\beta}}=\gamma_{\overline{\beta}}(I_{j(\beta)})$, where $I_{j(\beta)}$ is chosen such that $x_\beta^* \in  I_{j(\beta)}$.
Because the "cylinder sets" $I_{j(\beta)}$ are pairwise disjoints, we get 
$$ \sum_{\beta \in \mathscr{W}_n^j,\ \mathscr{D}(\alpha,\beta)<\epsilon} \Vert \gamma'_\beta \Vert_{I_j,\infty}^\delta
\leq C' \mu_\delta \left(\bigcup_{\beta \in \mathscr{W}_n^j,\ \mathscr{D}(\alpha,\beta)<\epsilon} I_{\overline{\beta}}\right).$$
We now conclude the proof by contemplating the implications of having $\mathscr{D}(\alpha,\beta)<\epsilon$. Roughly speaking, it implies that the {\it repelling fixed points} of the maps $\gamma_\alpha$ and $\gamma_\beta$ are $\epsilon$-close. Indeed, since we have
$$\mathscr{D}(\alpha,\beta)=2\inf_{x\in I_j} \frac{\vert c_\alpha d_\beta-c_\beta d_\alpha \vert}
{\vert c_\alpha x+d_\alpha\vert \vert c_\beta x +d_\beta \vert},$$
and
$$\gamma'_\alpha(x)=\frac{1}{(c_\alpha x+d_\alpha)^2},$$
we can use the bounded distortion property combined with (\ref{ineq1}) to observe that
$$ \mathscr{D}(\alpha,\beta)\geq \frac{1}{L} \left \vert \frac{d_\alpha}{c_\alpha}- \frac{d_\alpha}{c_\alpha}  \right \vert,$$
for some large constant $L>0$. Using (\ref{ineq2}) we deduce that 
$$\vert x_\alpha^*-x_\beta^* \vert \leq L' \left( \epsilon+e^{-\eta_0 n }\right).$$
Using the uniform contraction estimate, we get that the union of cylinder sets
$$ \bigcup_{\beta \in \mathscr{W}_n^j,\ \mathscr{D}(\alpha,\beta)<\epsilon} I_{\overline{\beta}}$$
is included in an interval of length at most $\widetilde{L}(\rho^n+\epsilon+e^{-\eta_0 n })$.
Choosing $\epsilon$ of size $e^{-\eta_1 n}$ with $\eta_1\leq \min\{\eta_0, \vert \log \rho \vert \}$ and using the estimate from 
proposition \ref{ADR} we conclude the proof. $\square$

\subsection{Proof of Proposition \ref{Decay1}}
We set $s=\sigma+it$, where $\sigma_0\leq \sigma\leq \delta$. We then write for $f \in C^1(I)$, 
$$\int_{\Lambda(\Gamma)} \vert \lt_{s,\theta}^n(f) \vert^2d\mu_\delta=
\sum_{j=1}^{2r} \int_{I_j} \sum_{\alpha,\beta \in \mathscr{W}_n^j} \left(\gamma'_\alpha\right)^{\sigma+it} \left (\gamma'_\beta\right)^{\sigma-it}
\chi_\theta(P\gamma_\alpha)\overline{\chi_\theta(P\gamma_\beta)} f\circ \gamma_\alpha \overline{f\circ \gamma_\beta} d\mu_\delta$$
$$=\sum_j \sum_{\alpha,\beta \in \mathscr{W}_n^j}\chi_\theta(P\gamma_\alpha)\overline{\chi_\theta(P\gamma_\beta)}
\int_{\Lambda(\Gamma)} e^{it\Phi_{\alpha,\beta}(x)} g_{\alpha,\beta}^{(j)}(x)d\mu_\delta(x), $$
where we have set
$$\Phi_{\alpha,\beta}(x):=\log\gamma'_\alpha(x)-\log \gamma'_\beta(x),$$
and
$$g_{\alpha,\beta}^{(j)}(x)=\varphi_j(x) \left(\gamma'_\alpha(x)\right)^{\sigma} \left (\gamma'_\beta(x)\right)^{\sigma}f\circ \gamma_\alpha(x) \overline{f\circ \gamma_\beta(x)},$$
with $\varphi_j$ being a $C^1(I)$ function which is $\equiv 1$ on $I_j$ and $\equiv 0$ on $I_i$ for $i\neq j$.

\bigskip
We point out that because they {\it do not depend} on the $x$ variable, but only on the word $\alpha$, the oscillating terms $\chi_\theta(P\gamma_\alpha)$ do not
interfere with the oscillatory integrals, which is the {\it crucial reason} why we will get estimates {\it uniform} with respect to $\theta$.

\bigskip
Using the bounded distortion estimate and the hyperbolicity estimate, we have
$$\Vert  g_{\alpha,\beta}^{(j)}\Vert_\infty \leq C_1 \Vert \gamma'_\alpha \Vert^\sigma_{\infty,j}\Vert \gamma'_\beta \Vert^\sigma_{\infty,j} \Vert f\Vert^2_{(t)},$$
while
$$\left \Vert  \frac{d}{dx}\left (g_{\alpha,\beta}^{(j)}\right) \right \Vert_\infty \leq C_2 \Vert \gamma'_\alpha \Vert^\sigma_{\infty,j}\Vert \gamma'_\beta \Vert^\sigma_{\infty,j} \Vert f\Vert^2_{(t)}(1+\vert t\vert \rho^n).$$
On the other hand we have precisely 
$$ \inf_{x\in I_j} \vert \Phi'_{\alpha,\beta}(x)\vert=\mathscr{D}(\alpha,\beta).$$
The bounded distortion estimates for the second (and third) derivatives show that 
$$\Vert \Phi_{\alpha,\beta} \Vert_{C^2(I_j)}\leq M$$
for some uniform $M>0$. We now pick $\epsilon=e^{-\eta n}$, with $0<\eta<\eta_0$ so that (UNI) holds and write
$$ \int_{\Lambda(\Gamma)} \vert \lt_{s,\theta}^n(f) \vert^2d\mu_\delta\leq \underbrace{
C_1 \sum_{j}\sum_{\mathscr{D}(\alpha,\beta)<\epsilon}\Vert \gamma'_\alpha \Vert^\sigma_{\infty,j}\Vert \gamma'_\beta \Vert^\sigma_{\infty,j} \Vert f\Vert^2_{(t)} 
}_{\mathrm{near\ diagonal\ sum}}+ \underbrace{\sum_{j}\sum_{\mathscr{D}(\alpha,\beta)\geq \epsilon} \left \vert  \int_{\Lambda(\Gamma)} e^{it\Phi_{\alpha,\beta}} g_{\alpha,\beta}^{(j)}d\mu_\delta  \right \vert}_{\mathrm{off\ diagonal\ sum}}.$$
Using the pressure estimate and the (UNI) property, the "near diagonal sum" is estimated from above by
$$C_3\Vert f\Vert^2_{(t)}A(\sigma_0,n)e^{nP(\sigma_0)}\epsilon^\delta.$$
Using the polynomial decay result on oscillatory integrals, the "off diagonal sum" is estimated from above (again using the pressure estimate) by
$$C_4\Vert f\Vert^2_{(t)}\frac{\vert t\vert^{-\beta_1}(1+\vert t\vert \rho^n)}{\epsilon^{\beta_2}}e^{2nP(\sigma_0)},$$
so that
$$\int_{\Lambda(\Gamma)} \vert \lt_{s,\theta}^n(f) \vert^2d\mu_\delta\leq C_5\Vert f\Vert^2_{(t)}
\left ( A(\sigma_0,n)e^{nP(\sigma_0)}\epsilon^\delta+   \frac{\vert t\vert^{-\beta_1}(1+\vert t\vert \rho^n)}{\epsilon^{\beta_2}}e^{2nP(\sigma_0)}\right).$$
We recall that $A(\sigma_0,n)\leq C\overline{\rho}^{(\sigma-\delta)n}$ and $\epsilon=e^{-\eta n}$. In the latter, $n$ is now taken as $n=[C_0\log\vert t\vert]$. We know \underline{fix} $C_0>>1$ so that $\vert t\vert \rho^n$ stays bounded as $\vert t\vert \rightarrow +\infty$ and choose
$\eta>0$ small enough so that we get
$$\int_{\Lambda(\Gamma)} \vert \lt_{s,\theta}^n(f) \vert^2d\mu_\delta\leq C_6\Vert f\Vert^2_{(t)}
\left ( A(\sigma_0,n)e^{nP(\sigma_0)}\vert t\vert^{-\beta_3}+   \vert t\vert^{-\beta_1 /2}e^{2nP(\sigma_0)}\right).$$
It is now clear that by taking $\sigma_0$ close enough to $\delta$ we obtain for all $\vert t \vert$ large,
$$ \int_{\Lambda(\Gamma)} \vert \lt_{s,\theta}^n(f) \vert^2d\mu_\delta\leq C_7\Vert f\Vert^2_{(t)}\vert t \vert^{-\beta_4},$$
for some $\beta_4>0$
and the proof is complete. $\square$

\section{Zeros of $\mathcal{Z}_\Gamma(s,\theta)$ on the line $\{ \Re(s)=\delta\}$}
In this final section we prove Proposition \ref{LowF}, by combining standard ideas from \cite{Naud2} and \cite{ParryPollicott2}. Notice that we already know from
\cite{Naud2}, that $s\mapsto \mathcal{Z}_\Gamma(s,0)$ only vanishes at $s=\delta$ on the line $\{ \Re(s)=\delta\}$, with a simple zero, a consequence of the fact that the length spectrum of $\Gamma \backslash \hh$ is non-lattice. Therefore we need to show
that on the line $\{ \Re(s)=\delta\}$, if $\theta\neq 0\ \mathrm{mod}\ \Z^r$, then
$s \mapsto \mathcal{Z}_\Gamma(s,\theta)$ does not vanish. Assume that $\theta\neq 0\ \mathrm{mod}\  \Z^r$ and suppose that 
$\mathcal{Z}_\Gamma(\delta+it_0,\theta)=0$. Then by the Fredholm determinant identity, we know that there exists $g=g_{t_0,\theta}\in C^1(I)$  with
$\Vert g\Vert_\infty\neq 0$ such that 
$$\lt_{\delta+it_0,\theta}(g)=g.$$
Using the Ruelle-Perron-Frobenius theorem, we can conjugate $\lt_{\delta+it_0,\theta}$ by the positive non vanishing eigenfunction $h_\delta$ so that we have
for all $x\in I_i$ 
$$\sum_{j\neq i}h_j(x)=1,$$
$$\sum_{j\neq i} h_j(x) (\gamma_j'(x))^{it}\chi_\theta(P\gamma_j) \widetilde{g}\circ \gamma_j(x)=\widetilde{g}(x),$$
where $h_j(x):=h_\delta^{-1} (\gamma'_j)^\delta h_\delta\circ \gamma_j$ and $\widetilde{g}=h_\delta^{-1}g$.
Choosing $i$ and $x_0 \in I_i$ such that 
$$\vert \widetilde{g}(x_0)\vert=
\sup_{x\in \Lambda(\Gamma)} \vert \widetilde{g}(x)\vert:=\Vert \widetilde{g} \Vert_{\infty,\Lambda(\Gamma)},$$
we observe that 
$$\vert \widetilde{g}(x_0)\vert \leq \sum_{j\neq i} h_j(x_0)\vert \widetilde{g}\circ \gamma_j(x_0)\vert \leq \Vert \widetilde{g}\Vert_{\infty,\Lambda(\Gamma)},$$
which implies that for all $j\neq i$, 
$$\vert \widetilde{g}\circ \gamma_j(x_0)\vert=\Vert \widetilde{g} \Vert_{\infty,\Lambda(\Gamma)}.$$
Iterating this argument and using the fact that the orbit of $x_0$ under the semigroup generated by $\gamma_1,\ldots,\gamma_{2r}$ is dense in $\Lambda(\Gamma)$,
we conclude that $\vert \widetilde{g}\vert$ is actually constant on $\Lambda(\Gamma)$. Taking this constant equal to one, we can write 
$$\widetilde{g}(x)=e^{i\phi(x)},$$ where $\phi$ is in say $C^0(\Lambda(\Gamma))$. We obtain that for all $x\in I_i\cap \Lambda(\Gamma)$, 
$$\sum_{j\neq i} h_j(x)e^{it_0\log\gamma'_j(x)+2i\pi \langle \theta,P\gamma_j\rangle+\phi\circ \gamma_j(x)}=e^{i\phi(x)},$$
which by strict convexity of the unit circle implies that for all $j$,
$$t_0 \log \gamma'_j(x)+\phi\circ \gamma_j(x)-\phi(x) \in 2\pi\Z- 2\pi \langle \theta,P\gamma_j\rangle.$$
If $t_0=0$ then this implies (by evaluating at attracting fixed points of each $\gamma_j$) that for all $j=1,\dots,r$,
$$\langle \theta,P\gamma_j\rangle\in \Z.$$
Since $\{P\gamma_1,\ldots,P\gamma_r \}$ is a $\Z$-basis of $\Z^r$, it implies that $\theta \in \Z^r$, a contradiction.
Therefore $t_0\neq 0$. Iterating the above formula, we get that for all $\gamma_\alpha \in \mathscr{W}^j_n$, 
$$t_0\log \gamma'_\alpha(x)+\phi\circ\gamma_\alpha(x)-\phi(x)\in 2\pi\Z-2\pi\langle \theta, P\gamma_\alpha\rangle.$$
By evaluating at the attracting fixed point $x_\alpha$ of $\gamma_\alpha$, we obtain that the translation length $l_\alpha$ of $\gamma_\alpha$, given by the formula
$$e^{-l_\alpha}=\gamma_\alpha'(x_\alpha),$$ 
satisfies
$$ l_\alpha\in \frac{2\pi}{t_0}\Z +\frac{2\pi}{t_0}\langle \theta, P\gamma_\alpha\rangle.$$
In term of lengths of closed geodesics, it shows in particular that the set of closed geodesics 
which belong to the homology class of $0$ (i.e. $P\gamma_\alpha=0$) is a subset of $\frac{2\pi}{t_0}\Z$. But this would imply that
the length spectrum of $(\mathrm{Ker}P) \backslash \hh$ is {\it lattice}, which is impossible since $\mathrm{Ker}P=[\Gamma,\Gamma]$ is the commutator subgroup of $\Gamma$ and hence non-elementary, see for example \cite{Dalbo}.

 \end{document}